\documentclass[leqno]{article}

\def\diam{\mathop{\rm diam}}
\def\dist{\mathop{\rm dist}}

\begin{document}

\title{Potpourri, 9}

\author{Stephen William Semmes	\\
	Rice University		\\
	Houston, Texas}

\date{}

\maketitle

	Let $(M, d(x, y))$ be a metric space.  Thus $M$ is a nonempty
set and $d(x, y)$ is a real-valued function defined for $x, y \in M$
such that $d(x, y) \ge 0$ for all $x, y \in M$, $d(x, y) = 0$ if and
only if $x = y$, $d(x, y) = d(y, x)$ for all $x, y \in M$, and
\begin{equation}
	d(x, z) \le d(x, y) + d(y, z)
\end{equation}
for all $x, y, z \in M$.

	A subset $E$ of $M$ is said to be \emph{bounded} if $d(x, y)$,
$x, y \in E$, is a bounded set of real numbers.  In this event we
define the \emph{diameter} of $E$, denoted $\diam E$, to be the
supremum of $d(x, y)$, $x, y \in E$.  This is interpreted as being $0$
if $E$ is the emptyset.

	If $E_1$, $E_2$ are subsets of $M$ with $E_1 \subseteq E_2$
and $E_2$ is bounded, then $E_1$ is bounded and
\begin{equation}
	\diam E_1 \le \diam E_2.
\end{equation}
If $E$ is a bounded subset of $M$ and $\overline{E}$ is the
\emph{closure} of $E$, i.e., the set of points in $M$ which are
elements of $E$ or limit points of $E$, then $\overline{E}$
is also a bounded subset of $M$ and
\begin{equation}
	\diam \overline{E} = \diam E.
\end{equation}

	Suppose that $(N, \rho(u, v))$ is another metric space, and
that $f$ is a mapping from $M$ to $N$.  We say that $f$ is
\emph{Lipschitz} if there is a nonnegative real number $C$ such that
\begin{equation}
	\rho(f(x), f(y)) \le C \, d(x, y)
\end{equation}
for all $x, y \in M$.  We may say that $f$ is $C$-Lipschitz in this
case to be explicit about the Lipschitz constant $C$.

	If $f(x)$ is a real-valued function on $M$, then $f$ is
$C$-Lipschitz if and only if
\begin{equation}
	f(x) \le f(y) + C \, d(x, y)
\end{equation}
for all $x, y \in M$.  In particular, if $p \in M$, then the function
\begin{equation}
	f_p(x) = d(x, p)
\end{equation}
is $1$-Lipschitz, by the triangle inequality for $d(z, w)$.

	More generally, suppose that $A$ is a nonempty subset of $M$,
and define the distance from a point $x \in M$ to $A$ by
\begin{equation}
	\dist (x, A) = \inf \{d(x, a) : a \in A\}.
\end{equation}
It is easy to see that
\begin{equation}
	\dist(x, A) \le \dist(y, A) + d(x, y)
\end{equation}
for all $x, y \in M$, so that the distance to $A$ defines a function
on $M$ which is $1$-Lipschitz.

	Let $E$ be a bounded subset of $M$, and let $f : M \to N$ be a
$C$-Lipschitz mapping from $M$ into $N$.  It is easy to see that
$f(E)$ is a bounded subset of $N$, and that
\begin{equation}
	\diam f(E) \le C \, \diam E,
\end{equation}
where of course the diameter of $f(E)$ uses the metric on $N$ while
the diameter of $E$ uses the metric on $M$.  This property is in fact
equivalent to saying that $f$ is $C$-Lipschitz.

	Let $E$ be a subset of $M$, and let $\alpha$ be a positive
real number.  We define $\mu^\alpha(E)$, the
$\alpha$-dimensional Hausdorff content of $E$, to be the infimum of
\begin{equation}
	\sum_{A \in \mathcal{A}} (\diam A)^\alpha
\end{equation}
over all families $\mathcal{A}$ of bounded subsets of $M$ such that
$\mathcal{A}$ has at most finitely many elements and
\begin{equation}
	E \subseteq \bigcup_{A \in \mathcal{A}} A,
\end{equation}
which is to say that $\mathcal{A}$ is a covering of $E$.

	If $p \in M$ and $r$ is a positive real number, the open and
closed balls in $M$ centered at $p$ and with radius $r$ are given by
\begin{equation}
	B(p, r) = \{x \in M : d(x, p) < r\}
\end{equation}
and
\begin{equation}
	\overline{B}(p, r) = \{x \in M : d(x, p) \le r\}.
\end{equation}
For any $p \in M$, $M$ is equal to the union of $B(p, l)$ as $l$ runs
through the set of positive integers, and thus every subset of $M$ is
contained in the union of a family of bounded subsets of $M$ with at
most countably many elements.

	It may be that $\mu^\alpha(E) = +\infty$.  We can also allow
unbounded sets in the covering, with the convention that the diameter
of an unbounded set is equal to $+\infty$.  If $E$ is a bounded subset
of $M$, then we automatically have that
\begin{equation}
	\mu^\alpha(E) \le (\diam E)^\alpha,
\end{equation}
just by using the covering of $E$ by $E$ itself.

	Notice that $\mu^\alpha(E) = 0$ if $E$ is the empty set, or if
$E$ contains just one element.  If $E$, $\widetilde{E}$ are subsets of
$M$ with
\begin{equation}
	\widetilde{E} \subseteq E,
\end{equation}
then
\begin{equation}
	\mu^\alpha(\widetilde{E}) \le \mu^\alpha(E)
\end{equation}
automatically.

	If $E_1$, $E_2$ are subsets of $M$, then it is easy to see
that
\begin{equation}
	\mu^\alpha(E_1 \cup E_2) 
		\le \mu^\alpha(E_1) + \mu^\alpha(E_2).
\end{equation}
Moreover, if $E_1, E_2, E_3, \ldots$ is any sequence of subsets of
$M$, then
\begin{equation}
	\mu^\alpha \bigg(\bigcup_{j=1}^\infty E_j \bigg)
		\le \sum_{j=1}^\infty \mu^\alpha(E_j),
\end{equation}
as one can verify.

	In the definition of $\mu^\alpha(E)$, one may as well
restrict one's attention to coverings $\mathcal{A}$ of $E$ by subsets
$A$ of $M$ which are closed.  Indeed, given any covering of $E$ by at
most countably many subsets of $M$, one can simply take the closures
of the subsets of $M$ in the covering to get a covering by closed
subsets of $M$.  The diameter of the closure of a subset of $M$ is
equal to the diameter of the original subset of $M$, and so the sum
employed in the definition of $\mu^\alpha(E)$ is not affected
by passing to the closures of the sets in the coverings.

	Let $A$ be any subset of $M$ and let $r$ be a positive real
number, and put
\begin{equation}
	A(r) = \bigcup \{B(x, r) : x \in A\}
		= \{x \in M : \dist(x, A) < r\}.
\end{equation}
This is an open subset of $M$ which contains $A$.

	If $A$ is a bounded subset of $M$, then it is easy to see
that $A(r)$ is a bounded subset of $M$, and that
\begin{equation}
	\diam A(r) \le \diam A + 2 \, r.
\end{equation}
Using this, one can show that one can restrict to coverings by open
subsets of $M$ in the definition of $\mu^\alpha(E)$ and get
the same result.

	Now suppose that $E$ is a compact subset of $M$.  In this
event one can restrict to finite coverings of $E$ in the definition of
$\mu^\alpha(E)$ and get the same result.  For one can first
reduce to the case of open coverings, as in the preceding paragraph,
and then to finite coverings by compactness.

	Let $\epsilon > 0$ be given.  By an \emph{$\epsilon$-family}
of subsets of $M$ we mean a family of subsets of $M$ which has at most
countably many elements and where each element of the family has
diameter less than $\epsilon$.

	Recall that a metric space $M$ is said to be \emph{separable}
if there is a subset of $M$ which is at most countable and also dense
in $M$.  This is equivalent to saying that for each $\epsilon > 0$
there is an $\epsilon$-family of subsets of $M$ such that the union of
the subsets of $M$ in the family is equal to all of $M$.

	Let $\alpha, \epsilon > 0$ be given, and let $E$ be a subset
of $M$.  Define $\mathcal{H}^\alpha_\epsilon(E)$ to be the infimum of
\begin{equation}
	\sum_{A \in \mathcal{A}} (\diam A)^\alpha
\end{equation}
over all $\epsilon$-families $\mathcal{A}$ of subsets of $M$ such that
\begin{equation}
	E \subseteq \bigcup_{A \in \mathcal{A}} A.
\end{equation}
If no such $\epsilon$-family exists, then put
$\mathcal{H}^\alpha_\epsilon(E) = +\infty$.

	If $0 < \epsilon_1 < \epsilon_2$, then
\begin{equation}
	\mathcal{H}^\alpha_{\epsilon_2}(E) 
		\le \mathcal{H}^\alpha_{\epsilon_1}(E).
\end{equation}
Basically $\mu^\alpha(E)$ corresponds to
$\mathcal{H}^\alpha_\epsilon(E)$ with $\epsilon = + \infty$, and in
particular we have that
\begin{equation}
	\mu^\alpha(E) \le \mathcal{H}^\alpha_\epsilon(E)
\end{equation}
for all $\epsilon > 0$.

	Recall that a subset $E$ of $M$ is said to be \emph{totally
bounded} if for each $\epsilon > 0$ there is a finite collection of
subsets of $M$ with diameter less than $\epsilon$ such that $E$ is
contained in the union of these subsets of $M$.  If $E$ is a totally
bounded subset of $M$, then $\mathcal{H}^\alpha_\epsilon(E)$ is finite
for all $\alpha, \epsilon > 0$.

	If $E$ is the empty set, or if $E$ contains just one element,
then $\mathcal{H}^\alpha_\epsilon(E) = 0$ for all $\alpha, \epsilon >
0$.  If $E$, $\widetilde{E}$ are subsets of $M$ with $\widetilde{E}
\subseteq E$, then
\begin{equation}
	\mathcal{H}^\alpha_\epsilon(\widetilde{E}) 
		\le \mathcal{H}^\alpha_\epsilon(E)
\end{equation}
for all $\alpha, \epsilon > 0$.  If $E_1, E_2, \ldots$ is a sequence
of subsets of $M$, then
\begin{equation}
	\mathcal{H}^\alpha_\epsilon \bigg(\sum_{j=1}^\infty E_j \bigg)
		\le \sum_{j=1}^\infty \mathcal{H}^\alpha_\epsilon(E_j)
\end{equation}
for all $\alpha, \epsilon > 0$.

	Just as for $\mu^\alpha$, one can restrict one's attention to
$\epsilon$-families of open or closed subsets of $M$ in the definition
of $\mathcal{H}^\alpha_\epsilon(E)$ and get the same result.  If $E$
is a compact subset of $M$, it follows that one can restrict one's
attention to $\epsilon$-families which contain only finitely many
subsets of $M$ in the definition of $\mathcal{H}^\alpha_\epsilon(E)$
and obtain the same result.

	Let $\epsilon > 0$ be given, and suppose that $E_1$, $E_2$ are
subsets of $M$ such that
\begin{equation}
	d(x, y) \ge \epsilon
\end{equation}
for all $x \in E_1$, $y \in E_2$.  In this event we actually have that
\begin{equation}
	\mathcal{H}^\alpha_\epsilon(E_1 \cup E_2)
		= \mathcal{H}^\alpha_\epsilon(E_1)
			+ \mathcal{H}^\alpha_\epsilon(E_2)
\end{equation}
for all $\alpha > 0$.  Indeed, if $A$ is any subset of $M$ with
diameter less than $\epsilon$, then $A$ may intersect either of $E_1$,
$E_2$ but not both.

	For each $\alpha > 0$ and $E \subseteq M$ define
$\mathcal{H}^\alpha(E)$ to be the supremum of
$\mathcal{H}^\alpha_\epsilon(E)$ over all $\epsilon > 0$, which is the
same as the limit as $\epsilon \to 0$.  This is called the
$\alpha$-dimensional Hausdorff measure of $E$, and it may be equal to
$+\infty$.

	A subset $E$ of $M$ satisfies $\mu^\alpha(E) = 0$ if and only
if for each $\eta > 0$ there is a family of subsets $\mathcal{A}$ of
$M$ with at most countably many elements such that
\begin{equation}
	E \subseteq \bigcup_{A \in \mathcal{A}} A
\end{equation}
and
\begin{equation}
	\sum_{A \in \mathcal{A}} (\diam A)^\alpha < \eta.
\end{equation}
In this event the sets in the families used to cover $E$ should have
small diameter, and we have that
\begin{equation}
	\mathcal{H}^\alpha_\epsilon(E) = 0
\end{equation}
for all $\epsilon > 0$, and hence
\begin{equation}
	\mathcal{H}^\alpha(E) = 0.
\end{equation}

	If $E$, $\widetilde{E}$ are subsets of $M$ with $\widetilde{E}
\subseteq E$, then we have
\begin{equation}
	\mathcal{H}^\alpha(\widetilde{E}) \le \mathcal{H}^\alpha(E)
\end{equation}
for all $\alpha > 0$.  Similarly, if $E_j$, $j \ge 1$, is a sequence
of subsets of $M$, then
\begin{equation}
	\mathcal{H}^\alpha\bigg(\bigcup_{j=1}^\infty E_j \bigg)
		\le \sum_{j=1}^\infty \mathcal{H}^\alpha(E_j)
\end{equation}
for all $j$, as a consequence of the analogous property for
$\mathcal{H}^\alpha_\epsilon$ for all $\epsilon > 0$.  If $E_1$, $E_2$
are subsets of $M$ and $r$ is a positive real number such that
\begin{equation}
	d(x, y) \ge r
\end{equation}
for all $x \in E_1$ and $y \in E_2$, then
\begin{equation}
	\mathcal{H}^\alpha (E_1 \cup E_2)
		= \mathcal{H}^\alpha(E_1) + \mathcal{H}^\alpha(E_2)
\end{equation}
for all $\alpha > 0$.

	If $0 < \alpha \le \beta$, $\epsilon > 0$, and $E$ is a subset
of $M$, then
\begin{equation}
	\mathcal{H}^\beta_\epsilon(E) 
	   \le \epsilon^{\beta - \alpha} \, \mathcal{H}^\alpha_\epsilon(E).
\end{equation}
In particular, if $\alpha < \beta$ and
\begin{equation}
	\mathcal{H}^\alpha(E) < \infty,
\end{equation}
then
\begin{equation}
	\mathcal{H}^\beta(E) = 0.
\end{equation}

	Suppose that $N$ is another metric space, and that $f$
is a $C$-Lipschitz mapping from $M$ to $N$ for some $C > 0$.
It is easy to see from the definitions that
\begin{equation}
	\mu^\alpha(f(E)) \le C^\alpha \, \mu^\alpha(E)
\end{equation}
for all $\alpha > 0$, where of course $\mu^\alpha(E)$ is defined using
the metric on $M$ and $\mu^\alpha(f(E))$ is defined using the metric
on $N$.
Similarly,
\begin{equation}
	\mathcal{H}^\alpha_{C \, \epsilon}(f(E))
		\le C^\alpha \, \mathcal{H}^\alpha_\epsilon(E)
\end{equation}
for all $\alpha, \epsilon > 0$, and hence
\begin{equation}
	\mathcal{H}^\alpha(f(E)) \le C^\alpha \, \mathcal{H}^\alpha(E).
\end{equation}

	Fix a $C$-Lipschitz real-valued function $f$ on $M$ and a
subset $E$ of $M$.  For each real number $t$, put
\begin{equation}
	E_t = \{x \in E : f(x) = t\}.
\end{equation}

	Let $\mathcal{A}$ be an at most countable family of nonempty
bounded subsets of $M$ such that
\begin{equation}
	E \subseteq \bigcup_{A \in \mathcal{A}} A.
\end{equation}
For each $t \in {\bf R}$, let $\mathcal{A}_t$ be the family of $A \in
\mathcal{A}$ such that
\begin{equation}
	A \cap E_t \ne \emptyset.
\end{equation}
Thus
\begin{equation}
	E_t \subseteq \bigcup_{A \in \mathcal{A}_t} A,
\end{equation}
where this is interpreted as automatic when $E_t = \emptyset$.

	For each $A \in \mathcal{A}$, $f(A)$ is a nonempty bounded
subset of the real line with diameter less than or equal to $C$ times
the diameter of $A$.  Let $I(A)$ denote the smallest closed interval
in the real line which contains $f(A)$, which is to say that $I(A) =
[a, b]$ with $a = \inf f(A)$ and $b = \sup f(A)$.

	Fix a real number $\alpha > 1$, and for each $t \in {\bf R}$
put
\begin{equation}
	h(t) = \sum_{A \in \mathcal{A}} (\diam A)^{\alpha - 1} \,
				\chi_{I(A)}(t).
\end{equation}
As usual, $\chi_{I(A)}(t)$ denotes the characteristic function of
$I(A)$ on the real line, which is equal to $1$ when $t \in I(A)$ and
is equal to $0$ otherwise.

	Thus for each $t \in {\bf R}$ we have that $h(t) \ge 0$ by
definition, and $h(t)$ may be equal to $+ \infty$.  If $\mathcal{A}$
contains only finitely many subsets of $M$, then $h(t)$ is a finite
step function on the real line.  If $E$ is a compact subset of $M$,
then it is natural to restrict our attention to finite families
$\mathcal{A}$ of subsets of $M$, as before.

	A basic feature of this function $h(t)$ is that
\begin{equation}
   \int_{\bf R} h(t) \, dt \le C \sum_{A \in \mathcal{A}} (\diam A)^\alpha.
\end{equation}
As usual, this is somewhat simpler when $\mathcal{A}$ is a finite
family, as when $E$ is compact.

	For each $t \in \mathcal{R}$ we have that
\begin{equation}
	\mu^{\alpha - 1}(E_t) \le h(t).
\end{equation}
If $\epsilon > 0$ and $\mathcal{A}$ is an $\epsilon$-family of subsets
of $M$, then
\begin{equation}
	\mathcal{H}^{\alpha - 1}_\epsilon(E_t) \le h(t)
\end{equation}
for all $t \in {\bf R}$.

	For instance, suppose that $\mu^\alpha(E) = 0$.  For each $r >
0$, consider the set of $t \in {\bf R}$ such that $\mu^{\alpha -
1}(E_t) \ge r$.  Let $\eta > 0$ be given, and choose $\mathcal{A}$ so
that $\sum_{A \in \mathcal{A}} (\diam A)^\alpha < \eta$.  Thus the
integral of $h$ over the real line is less than $C \, \eta$.  One can
use this to show that the set of $t \in {\bf R}$ such that
$\mu^{\alpha - 1}(E_t) > 0$ has measure $0$ in the real line.

	More generally, suppose that $\mu^\alpha(E) < \infty$.  For
each positive integer $j$, let $\mathcal{A}_j$ be an at most countable
family of subsets of $M$ such that $E \subseteq \bigcup_{A \in
\mathcal{A}_j} A$ and $\sum_{A \in \mathcal{A}_j} (\diam A)^\alpha$ is
less than $\mu^\alpha(E) + (1/j)$.  For each $j$, also let $h_j(t)$ be
the function on the real line as in the preceding discussion.  If we
put $\phi(t) = \inf \{h_j(t) : j \ge \}$, then $\mu^{\alpha - 1}(E_t)
\le \phi(t)$ for all $t \in {\bf R}$.  Furthermore, the integral of
$\phi(t)$ over the real line is less than or equal to $\mu^\alpha(E)$.

	Now suppose that $\mathcal{H}^\alpha(E) < \infty$.  For each
$j$, let $\mathcal{A}$ be an $\epsilon$-family of subsets of $M$ with
$\epsilon < 1/j$ such that $E \subseteq \bigcup_{A \in \mathcal{A}_j}
A$ and $\sum_{A \in \mathcal{A}_j} (\diam A)^\alpha$ is less than
$\mathcal{H}^\alpha(E)$ plus $1/j$.  For each $j$ let $h_j(t)$ be the
function on the real line associated to $\mathcal{A}_j$ as in the
earlier discussion.  If $\psi(t)$ is equal to $\liminf_{j \to \infty}
h_j(t)$ for each real number $t$, then we have that
$\mathcal{H}^{\alpha - 1}(E_t)$ is less than or equal to $\psi(t)$ for
all $t \in {\bf R}$.  Also, the integral of $\psi(t)$ over the real
line is less than or equal to $\mathcal{H}^\alpha(E)$, by Fatou's
lemma.

\end{document}